\documentclass{svproc}
\usepackage{url}
\usepackage{amsthm}

\usepackage{amsmath,amssymb,amsfonts}
\usepackage{algorithmic}
\usepackage{graphicx}
\usepackage{textcomp}
\usepackage{xcolor}
\usepackage[misc]{ifsym}

\newtheorem{theo}{Theorem}
\newtheorem{lem}{Lemma}
\theoremstyle{remark}
\newtheorem{rmk}{Remark}
\newtheorem{exm}{Example}

\begin{document}
\mainmatter              
\title{Decentralized design of leader-following consensus protocols for asymmetric matrix-weighted heterogeneous multiagent systems}
\titlerunning{Decentralized design of consensus protocols}  
%
\author{Lanhao Zhao\inst{1,2} \and Yangzhou Chen\inst{1,2}$^{(\textrm{\Letter})}$
}
\authorrunning{Lanhao Zhao et al.} 
%
%
\institute{College of Artificial Intelligence and Automation, Beijing University of Technology.\\
\and
Engineering Research Center of Digital Community, Ministry of Education.\\
\email{zhaolanhao@emails.bjut.edu.cn, yzchen@bjut.edu.cn (Corresponding author)}}

\maketitle              

\begin{abstract}
This paper investigates a decentralized design approach of leader-following consensus protocols for heterogeneous multiagent systems under a fixed communication topology with a directed spanning tree (DST) and asymmetric weight matrix. First, a control protocol using only the information of the neighbor on the DST of each agent is designed, which is called the consensus protocol with minimal communication links. Particularly, the DST-based linear transformation method is used to transform the consensus problem into a partial variable stability problem of a corresponding system, and a decentralized design method is proposed to find the gain matrices in the protocols. Next, the decentralized design approach is extended to the protocols using all neighbor information in the original communication topology with the help of the matrix diagonally dominant method. Some numerical simulations are given to illustrate the theoretical results.
\end{abstract}
\section{Introduction}
The consensus of multi-agent systems (MASs) has received extensive attention in the control community in the past three decades (see, e.g., the review paper \cite{1}).
However, the following challenging issues in consensus problems are still necessary to further investigate: decentralized design approach of consensus protocols for heterogeneous MASs and generic asymmetric matrix-weighted communication topology.

In terms of consensus protocols for heterogeneous MASs, two key problem need further attention: one is the decentralized design approach of the feedback gains based on the distributed structure of consensus protocols. The other is how to use as little communication as possible from the neighbors of each agent. For heterogeneous MAS, the consensus problem of heterogeneous MAS with low-order individuals was first discussed and many results were obtained (see, e.g., \cite{12}).
For the heterogeneous MAS with high-order individuals, Wahrburg and Adamy \cite{w1} obtained sufficient and necessary conditions for state consensus problems by compensating for parameter deviations to isomorphize agents.
However, isomorphism has strict premises and is not always feasible.
Tian and Zhang \cite{t1} discussed the consensus problem of high-order heterogeneous MASs with unknown communication delay based on a class of transfer function models.
However, some severer assumptions are made in \cite{t1}, for example, the MAS is required to be semi-stable. 
Under more relaxed assumptions, our previous works for leaderless \cite{c1} and leader-following \cite{31} proposed a new method to deal with the consensus problem of heterogeneous MASs.
Based on a linear transformation constructed according to the characteristics of MAS, the consensus problem of heterogeneous MAS is transformed into the partial variable stability problem of the corresponding system \cite{34}, and some necessary and sufficient conditions were obtained.
It should be noted that the method developed in \cite{31} allowed that each agent has its own exclusive feedback gain. However, global information is still required when finding these feedback gains, which is difficult for large-scale heterogeneous MASs. 
To address this problem, a fully distributed method was proposed in \cite{f} to deal with this problem, where although global information is not required, but the feedback gain of all agents was set to be the same.
A natural extension is to allow each individual to have an exclusive feedback gain matrix.
This will increase the degree of freedom in the design of the gain matrix, and it is especially needed for heterogeneous MASs.
For the second problem mentioned above, researchers usually use intermittent communication methods such as sampled control \cite{32} and event-triggered control \cite{33} to reduce communication traffic.
On the other hand, using information from some neighbors rather than all neighbors is also a method of reducing communication \cite{ccdc}.
Thus, this prompts us to consider a consensus protocol that only uses part of neighbor information and design feedback gain matrices decentrally.

In addition, in practical applications, the weights expressing communication relationships between different agents may have different values concerning different components. 
At this point, weight in matrix form is needed to describe the interaction between two agents.
The consensus problem with communication weights in matrix form has been addressed in our early research \cite{24}, although this topic was not particularly emphasized there.
The concept of the matrix-weighted network has also been proposed by Trinh et al.\cite{23} when studying bearing-based formation control problems.
Recently, the consensus problem with the matrix-weighted network has received further attention, and some significant results have been obtained \cite{27}.
However, in most existing research, the weighted matrices are assumed to be symmetric \cite{28}, either positive definite or negative definite.
Compared with a symmetric weighted matrix, the asymmetric weighted matrix can describe a wider range of networks.
Asymmetric weighted matrices will cause difficulty in stability analysis in the work mentioned above, which is another key issue that inspires the research in this paper.
Therefore, a natural further consideration is to handle the consensus problem of heterogeneous MASs with asymmetric weighted matrices.

In this paper, we attempt to propose a decentralized design approach of leader-following consensus protocols for heterogeneous MASs with asymmetric weighted matrices.
The main contributions are summarized below.

1) For the leader-following consensus problem of heterogeneous MASs with asymmetric weighted matrices, the DST-based linear transformation method is used in the analysis and design of the control protocol, and the partial variable stability theory is adopted to obtain consensus criteria. 
Especially, the results obtained relax the assumption of symmetric weighted matrices in the existing work on matrix weighted MASs \cite{23,27,28}. 

2) A new consensus protocol is proposed, where each agent uses only the information of itself and its neighbors on the DST, and the gain matrices can be designed via a decentralized procedure.
This design approach is significant in reducing communication load and is called the consensus protocol design with minimal communication links.

3) The decentralized design approach is extended to find the gain matrices of the consensus protocol using all neighbors' information by adopting a matrix diagonally dominant idea.
Compared to previous work \cite{31}, the decentralized design approach provides computational convenience in designing feedback gains.

As far as we know, this is the first time to give the decentralized design method of consensus protocols for heterogeneous interconnected MASs with asymmetric weighted matrices. 

The rest of this paper is organized as follows.
Section \uppercase\expandafter{\romannumeral2} presents preliminary knowledge and the model of heterogeneous MASs with asymmetric weighted matrices.
Section \uppercase\expandafter{\romannumeral3} shows the decentralized design method of the consensus protocol for heterogeneous MASs with asymmetric weighted matrices.
Section \uppercase\expandafter{\romannumeral4} provides some numerical examples.
The summary is made in Section \uppercase\expandafter{\romannumeral5}.
\section{Heterogeneous MASs with asymmetric weight matrix}
In this section, we introduce some basic concepts and the model of heterogeneous MAS with asymmetric weighted matrices.
Let $\mathbb{R}^n$ be the $n$-dimensional Euclidean space
and ${\bf e}_1=\left[\begin{array}{llll}1 & 0 & \cdots & 0\end{array}\right]^T \in \mathbb{R}^{N+1}$ be the unit vector. The symbol $\otimes$ is used to express the Kronecker product.

Consider a linear parameter heterogeneous MAS composed of  $N+1$ agents, where agent $0 $ is the leader and the other agents are the followers. The dynamic of each follower is
\begin{equation}
	\dot{x}_i=A_i x_i+B_i u_i, i = 1, \cdots, N
\end{equation}
and the dynamic of the leader is
\begin{equation}
	\dot{x}_0=A_0 x_0+B_0 u_0
\end{equation}
where for $i = 0,1, \cdots, N$, $A_i \in \mathbb{R}^{n \times n}, B_i \in \mathbb{R}^{n \times m}$, $x_i \in \mathbb{R}^n$, $u_i \in \mathbb{R}^m$. Here $u_i, i=1,\cdots, N$ are the control inputs of the followers to be designed and $u_0$ is the external input of the leader.

The communication topology of the MAS can be represented by a directed graph $\mathcal{G}=(V, E, W)$ with weight matrix $W$, where $V=\{0,1,2, \cdots, N\}, E \subseteq\{(j, i): i, j \in V\}$. The edge $(j, i) \in E$ indicates that the information can be transferred to the agent $i$ from the agent $j$, and thus the agent $j$ is called the neighbor of the agent $i$.
The set $N_i$ is used to represent the collection of all neighbors of agent $i$.
There is no communication edge from the followers to the leader, which means that the leader's state is not affected by the followers' states, but only by external inputs.
The weight matrix $W$ is expressed by the block adjacency matrix $W=\left[W_{i j}\right]_{i,j=0}^N \in \mathbb{R}^{n(N+1) \times n(N+1)}$, where $W_{i j} \in  \mathbb{R}^{n \times n}$ denotes the communication weight matrix between $i$ and $j$ and $W_{i j} \neq 0$ if and only if $(j, i) \in E$.

We emphasize that for the matrix-weighted MASs, the communication weights are allowed to be in the form of matrices with generic entries.
For example, the communication weights between different components of an individual and the corresponding components of its neighbors may be different.

For the directed graph  $\mathcal{G}$, a sequence of end-to-end directed edges in the same direction is called a directed path (strong path).
If there is a node $i$ that has a directed path to any other node $j \neq i$, the $\mathcal{G}$ is called containing a directed spanning tree (DST) with the root node $i$ and thus quasi-strongly connected.
Furthermore, if there exists a directed path from any node to any other node, then graph $\mathcal{G}$ is called strongly connected. 

Furthermore, to provide a clearer description of the connection between the leader and the followers, we introduce the block matrix $\delta=\left[D_1^T, \ldots, D_N^T\right]^T$, where $D_i=W_{i0} \in \mathbb{R}^{n \times n}$ is a non-zero matrix iff the leader can communicate directly with the follower $i$, otherwise $D_i=0$. 

In this paper, it is assumed that there is a DST $\mathcal{T}=\left(V, E_\mathcal{T}, W_\mathcal{T}\right)$ with $0$ as the root node. Each non-root node $i \in\{1,2, \cdots, N\}$ has only one parent node denoted by $k_i$ and, without loss of generality, we assume that $k_i < i$. 

\section{Decentralized design of consensus protocol}

In this section, we first discuss the decentralized design of the consensus protocol, where each agent only uses the information of itself and its neighbor on the DST and call it the consensus protocol with minimal communication links.
Then we extend the approach to the protocol that uses all neighbor information by means of the matrix diagonally dominant idea.  

\subsection{Decentralized design of consensus protocol with minimum communication links}

The well-known fact is that in consensus problems with fixed communication topologies,
the existence of a DST is a necessary communication condition for the agents to reach consensus.
From the perspective of graph connectivity requirements, DST generally has the smallest number of edges compared to the overall network.
Therefore, at least two benefits would be gained if each agent only used the information of its neighbor on the DST instead of using all neighbor information in the consensus protocol design.
The first one is to reduce the communication load by using as few communication connections as possible.
The other one is that, as one sees later, one can directly design the gain matrices by a decentralized approach. 

Therefore, we propose the consensus protocol as follows, where the follower $i$, $i=1, \cdots, N$, only uses its own state (if available) and the relative state with its neighbor on the DST
\begin{equation}
	\begin{aligned}
		u_i=-G_i x_i+K_i W_{i, k_i}\left(x_{k_i}-x_i\right)\\
	\end{aligned}
\end{equation}
and $G_i, K_i \in \mathbb{R}^{m \times n}$ are the gain matrices that need to be designed. Here $W_{i,k_i}=D_i \in \mathbb{R}^{n \times n}$ if $k_i=0$.

The consensus protocol (3) consists of two parts. The first part is the agent's own state feedback, where the gain matrix $G_i$ is directly set to a zero matrix if the agent's own state is not available.
The second part is the feedback of the relative state describing the effect of the neighbor information. 
\begin{rmk}
	The protocol (3) is referred to as the protocol with minimal communication links because it only uses neighbor information on the DST.
	At this point, the communication topology can be replaced by the DST, which is a necessary condition for the agents to reach consensus.
	If the number of communication connections used in the protocol is further reduced, it will inevitably damage the DST and may lead to the consensus not being reached. 
\end{rmk}

It follows from (1) and (3) that
\begin{equation}
	\begin{gathered}
		\dot{x}_i=\left(A_i-B_i G_i\right) x_i+B_i K_i W_{i, k_i}\left(x_{k_i}-x_i\right),
		i=1,\cdots,N
	\end{gathered}
\end{equation}
By combining (2) and (4), we rewrite the MAS as follows
\begin{equation}
	\left[\begin{array}{c}
		\dot{x}_0 \\
		\dot{x}
	\end{array}\right]=\left[\begin{array}{cc}
		A_0 & 0 \\
		B_{\mathrm{D}} K_{\mathrm{D}} \delta & M
	\end{array}\right]\left[\begin{array}{c}
		x_0 \\
		x
	\end{array}\right]+\tilde{B} u_0
\end{equation}
where $x=\left[x_1^T, \ldots, x_N^T\right]^T$,
$$
\begin{gathered}
	M=A_{\mathrm{D}}-B_{\mathrm{D}}\left(G_{\mathrm{D}}+K_{\mathrm{D}} L_{W_\mathcal{T}}+K_{\mathrm{D}} \Delta\right),
	A_{\mathrm{D}}=\operatorname{diag}\left(A_1, \cdots, A_N\right)\\ B_{\mathrm{D}}=\operatorname{diag}\left(B_1, \cdots, B_N\right),
	G_{\mathrm{D}}=\operatorname{diag}\left(G_1, \ldots, G_N\right), K_{\mathrm{D}}=\operatorname{diag}\left(K_1, \ldots, K_N\right) \\
	\Delta =\operatorname{diag}\left(D_1, \ldots, D_N\right), \tilde{B}=\left[B_0^T, 0_{n \times m}, \ldots, 0_{n \times m}\right]^T
\end{gathered}
$$
and $L_{W_\mathcal{T}}=\left[L_{i j}\right]_{i, j=1}^N$ is the block Laplacian concerning the followers' network in the DST and defined as
$$
L_{i j}= \begin{cases}W_{i, k_i}, & j=i \\ -W_{i, k_i}, & j \neq i, j = k_i \\ 0, & j \neq i, j \neq k_i\end{cases}
$$
According to the definition of $L_{W_\mathcal{T}}$ and the assumption $k_i<i$, we easily conclude that matrix $L_{W_\mathcal{T}}$ is lower triangular.

Therefore, the leader-following state consensus problem can be described as follows: to find the gain matrices $G_i$ and $K_i$ in the protocol (3) such that the closed-loop system (5) meets
\begin{equation}
	\begin{gathered}
		\lim _{t \rightarrow \infty}\left\|x_i(t)-x_0(t)\right\|=0, \;
		\forall i \in\{1, \cdots, N\}
	\end{gathered}
\end{equation}
In this case, it is said that the leader-following state consensus of MAS (1) and (2) with protocol (3) is achieved.

The decentralized design refers to finding the gain matrices $G_i$ and $K_i$ without using global information such as Laplacian matrices.

Next, we transform the consensus problem expressed by (5) into the partial variable stability problem of a corresponding system based on the DST-based linear transformation method \cite{36} and provide a decentralized design approach for the gain matrices.

Consider the incidence matrix $P_0$ of the DST $\mathcal{T}=\left(V, E_\mathcal{T}, W_\mathcal{T}\right)$  with the leader as the root node.
Let the edges of the DST be represented by $\left(k_i, i\right) \in E, i=1,2, \cdots, N$.
Let $e_{k_i, i} \in \mathbb{R}^{N+1}$ be the incidence vector of edge $\left(k_i, i\right)$, where $k_i$-th component is $1$, $i$-th component is $-1$, and the remaining components are $0$.
Thus, we have $P_0=\left[e_{k_1, 1}, e_{k_2, 2}, \cdots, e_{k_N, N}\right]$.

Construct the following transformation matrix
$$
P=\left[\begin{array}{ll}
	P_0 & \mathbf{e}_1
\end{array}\right]^T \otimes I_n=:\left[\begin{array}{cc}
	p & \tilde{P}_0 \\
	1 & 0
\end{array}\right] \otimes I_n
$$
and the corresponding linear transformation  
\begin{equation}
	\left[\begin{array}{c}
		y \\
		x_0
	\end{array}\right]=P\left[\begin{array}{c}
		x_0 \\
		x
	\end{array}\right]
\end{equation}

The linear transformation (7) turns system (5) into
\begin{equation}
	\begin{aligned}
		& {\left[\begin{array}{c}
				\dot{y} \\
				\dot{x}_0
			\end{array}\right]=\left(\left[\begin{array}{cc}
				p & \tilde{P}_0 \\
				1 & 0
			\end{array}\right] \otimes I_n\right) \cdot\left[\begin{array}{cc}
				A_0 & 0 \\
				B_{\mathrm{D}} K_{\mathrm{D}} \delta & M
			\end{array}\right]} 
		 \left(\left[\begin{array}{cc}
			0 & 1 \\
			\tilde{P}_0^{-1} & 1_N
		\end{array}\right] \otimes I_n\right)\left[\begin{array}{c}
			y \\
			x_0
		\end{array}\right]+\left[\begin{array}{c}
			\hat{B} \\
			B_0
		\end{array}\right] u_0
	\end{aligned}
\end{equation}
and thus we get
\begin{equation}
	\begin{aligned}
		& \dot{y}=\bar{A} y+\hat{A} x_0+\hat{B} u_0 
		, \dot{x}_0=A_0 x_0+B_0 u_0
	\end{aligned}
\end{equation}
where
$$
\begin{aligned}
	\bar{A}= & \left(\tilde{P}_0 \otimes I_n\right) M\left(\tilde{P}_0^{-1} \otimes I_n\right) \\
	\hat{A}= & {\left[\left(A_0-A_1+B_1 G_1\right)^T,\left(A_1-B_1 G_1-A_2+B_2 G_2\right)^T\right.} \\
	& \left.\ldots,\left(A_{N-1}-B_{N-1} G_{N-1}-A_N+B_N G_N\right)^T\right]^T \\
	\hat{B}= & {\left[B_0^T, 0_{n \times m}, \ldots,0_{n \times m}\right]^T }
\end{aligned}
$$
\begin{rmk}
	Both matrices $\bar{A}$ and $M$ have the same Hurwitz stability because $\bar{A}$ is obtained by a similar transformation of $M$.
\end{rmk}

Introducing the following notations
$$
\eta=\left[x_0^T, u_0^T\right]^T, \quad \bar{B}=\left[\begin{array}{ll}
	\hat{A} & \hat{B}
\end{array}\right], \quad \bar{D}=\left[\begin{array}{cc}
	A_0 & B_0 \\
	0 & 0
\end{array}\right]
$$
where  $\bar{B}$ is $N n \times(n+m)$ matrix and $\bar{D}$ is  $(n+m) \times(n+m)$ matrix,
the equation (9) becomes
\begin{equation}
	\begin{aligned}
		& \dot{y}=\bar{A} y+\bar{B} \eta \\
		& \dot{\eta}=\bar{D} \eta
	\end{aligned}
\end{equation}

We can get from (7) that $y=0$ if and only if $x_ 0 = x_ 1=\ldots=x_ N $. Therefore, we have the following lemma.

\begin{lem}
	The leader-following state consensus of MAS (5) can be achieved iff the zero equilibrium point of system (10) is globally asymptotically $y$-stable.
\end{lem}

Now, the consensus problem has been transformed into a partial variable stability problem of the corresponding system (10), and then we can give the $y$-stability condition of (10) based on the partial stability theory \cite{34}.

First, construct the observability matrix of $(\bar{B}, \bar{D})$ as follows
$$
V_k=\left[\begin{array}{ll}
	\bar{B}^T \quad(\bar{B} \bar{D})^T \quad \cdots\left(\bar{B} \bar{D}^k\right)^T
\end{array}\right]^T, \; k=0, \ldots,(n+m-1)
$$

Next, matrices $L_1 \in {\mathbb R}^{h\times (n+m)}$ and $L_3 \in {\mathbb R}^{(n+m) \times h}$ are obtained via the following steps.

1) Let $s=\min \{k: \operatorname{rank}\left.V_k=\operatorname{rank} V_{k+1}\right\}$ and $h=\operatorname{rank} V_s$. 

2) For matrix $V_s$, remain the first row if it is not zero, otherwise it will be removed. Then starting from the second row, remove all rows linearly related to the previously remained rows. All the remained rows form a new matrix, denoted as $L_1$.

3) The set of linearly independent columns of $L_1$, say,  the columns $i_1, \cdots, i_h$ of $L_1$, constitutes a new reversible matrix $L_2$.

4) The row $i_j$ of $L_3$ is the row $j$ of the inverse matrix of $L_2, j=1, \ldots, h$, and the rest rows of $L_3$ are set to be zero.

\begin{lem}
	The zero equilibrium state of the system (10) is globally asymptotically $y$-stable iff the following auxiliary system (11) is asymptotically stable
	\begin{equation}
		\dot{\xi}=\bar{M} \xi, \bar{M}=\left[\begin{array}{cc}
			\bar{A} & \bar{B} L_3 \\
			0 & L_1 \bar{D} L_3
		\end{array}\right]
	\end{equation}
\end{lem}

Lemma 1 reveals the relationship between the consensus problem of the MAS and the partial variable stability problem of a corresponding system. Lemma 2 provides a criterion of partial variable stability based on the auxiliary system.
Therefore, it follows from the two lemmas that the consensus of MAS (1) and (2) under control protocol (3) is reached if and only if $\bar{M}$ is Hurwitz.
We note that $\bar{M}$ is an upper triangular block matrix, and thus its Hurwitz stability is equivalent to that of $\bar{A}$ and $L_ 1 \bar {D} L_3 $.
It implies that the state consensus of MAS (1) and (2) under control protocol (3) is achieved if and only if the matrix $M$ in (5) and $L_ 1 \bar {D} L_ 3$ in (11) are both Hurwitz stable.

Let  ${A_i}^*=A_i-B_i G_i$ and ${B_i}^*=B_i K_i\left(W_{i k_i}+D_i\right)$.
Then  according to the characteristic that matrix $L_{W_{\cal T}}$ is a lower triangular matrix, the matrix $M$ can be written as a lower triangular block matrix
$$
M=\left[\begin{array}{cclc}
	{A_1}^*-{B_1}^* & 0 & \cdots & 0 \\
	* & {A_2}^*-{B_2}^* & \cdots & 0 \\
	\cdots & \cdots & \cdots & \cdots \\
	* & * & \cdots & {A_N}^*-{B_N}^*
\end{array}\right]
$$
and thus we get the following theorem.

\begin{theo}
	The leader-following state consensus of MAS (1) and (2) under control protocol (3) is achieved iff matrices ${A_i}^*-{B_i}^*, i = 1,\cdots,N$ and $L_1 \bar{D} L_3$ are Hurwitz stable. 
\end{theo}

Furthermore, we point out that the gain matrices $G_i$ and $K_i$ to be designed in protocol (3) have different roles.
First, $L_1 \bar{D} L_3$ expresses the heterogeneous characteristics of MAS (1) and (2) and only depends on the choice of gain matrix $G_i, i = 1,\cdots,N$.
So one can first choose $G_i, i = 1,\cdots,N$ (if state $x_i$ available) such that $L_1 \bar{D} L_3$ is Hurwitz stable.
Moreover, as shown in Example 2 in Section \uppercase\expandafter{\romannumeral5}, the gain matrix $G_i$ is related to the convergence rate of the MAS, so one can first select them according to the actual demand for the convergence rate.
Next, after $G_i, i = 1,\cdots,N$ being selected, one further determines $K_i$ such that ${A_i}^*-{B_i}^*$ for $i = 1,\cdots,N$ are Hurwitz stable.

\begin{rmk}
	It is obvious that the design of the gain matrix $K_i$ according to Theorem 1 only depends on the parameters of agent $i$ and thus the design procedure is decentralized. 
\end{rmk}

\subsection{Decentralized design of consensus protocol that uses all neighbor information}
In this subsection, we extend the decentralized design approach to the case of the protocol that uses all neighbor information using the matrix diagonally dominant idea.

Consider the following protocol that uses all the neighbors' information

\begin{equation}
	\begin{aligned}
		u_i=-G_i x_i+K_i \sum_{j \in N_i} W_{i,j}\left(x_j-x_i\right)\\
	\end{aligned}
\end{equation}

In this case, the block Laplacian $L_W=\left[L_{i j}\right]_{i, j=1}^N$ concerning the followers' network in the original communication topology is defined as
$$
L_{i j}=\left\{\begin{array}{lr}
	\sum_{k=1, k \neq i}^N W_{i k}, & j=i \\
	-W_{i j}, & j \neq i, j \in N_i \\
	0, & j \neq i, j \notin N_i
\end{array}\right.
$$

Similarly, we use the DST-based linear transformation to transform the consensus problem of MAS (1)  and (2) with the protocol (12) into the partial variable stability problem of the following system
\begin{equation}
	\begin{aligned}
		& \dot{y}=\bar{A}' y+\hat{A} x_0+\hat{B} u_0 \\
		& \dot{x}_0=A_0 x_0+B_0 u_0 
	\end{aligned}
\end{equation}
where
$$
\bar{A}'=\left(\tilde{P}_0 \otimes I_n\right) M'\left(\tilde{P}_0^{-1} \otimes I_n\right),
M'=A_{\mathrm{D}}-B_{\mathrm{D}}\left(G_{\mathrm{D}}+K_{\mathrm{D}} L_{W}+K_{\mathrm{D}} \Delta\right)
$$
We have the following similar result with the same matrices $\bar{B}, \bar{D}, L_1$, and $ L_3$ as before.
\begin{lem}
	The zero equilibrium state of the system (13) is globally asymptotically $y$-stable iff the following auxiliary system (14) is asymptotically stable:
	\begin{equation}
		\dot{\xi}=\bar{M}' \xi, \;\; \bar{M}'=\left[\begin{array}{cc}
			\bar{A}' & \bar{B} L_3 \\
			0 & L_1 \bar{D} L_3
		\end{array}\right]
	\end{equation}
\end{lem}

According to Lemma 3, the leader-following state consensus of MAS (1) and (2) under control protocol (12) is achieved iff $\bar{M}'$ is Hurwitz stable, or equivalently, if and only if  $\bar{A}'$ and $L_1 \bar{D} L_3$ are Hurwitz stable.
Furthermore, the two matrices $\bar{A}'$ and $M'$ have the same Hurwitz stability because $\bar{A}'$ is obtained by a similar transformation of $M'$.

Let $A'_i=A_i-B_i G_i, B'_i=B_i K_i\left(\sum_{k=1, k \neq i}^N W_{i k}+D_i\right)$. 
Then the matrix $M'$
can be written in the following form
$$
\left[\begin{array}{cccc}
	A'_1-B'_1 &-B_1 K_1 W_{1 2} &\cdots& -B_1 K_1 W_{1 N} \\
	-B_2 K_2 W_{2 1} & A'_2-B'_2&\cdots & -B_2 K_2 W_{2 N} \\
	\cdots & \cdots & \cdots & \cdots \\
	-B_N K_N W_{N1} & -B_N K_N W_{N2} & \cdots &A'_N-B'_N
\end{array}\right]
$$
Therefore, we have the following result using the Gerschgorlin circle theorem of block matrices \cite{35}, where the matrix norms are defined as 
$$
\|H\|=\max _{1 \leq i \leq m} \sum_{j=1}^n\left|h_{i j}\right|,\left(\left\|H^{-1}\right\|\right)^{-1}=\min _{1 \leq i \leq m} \sum_{j=1}^n\left|h_{i j}\right|
$$
for any matrix $H =\left[h_{i j}\right] \in \mathbb{C}^{m \times n}$.

\begin{theo}
	The leader-following state consensus of MAS (1) and (2) with control protocol (12) is achieved if matrix $L_1 \bar{D} L_3$ is Hurwitz and the gain matrices $K_i,i=1,..., N $ in (12) are selected such that the following Gerschgorlin circles
	\begin{equation}
		\left(\left\|\left(A_i^{\prime}-B_i^{\prime}-\lambda I_n\right)^{-1}\right\|\right)^{-1} \leqq \sum_{ j \neq i, j \in N_i}\left\|-B_i  K_i W_{i,j}\right\|
	\end{equation}
	are all located in the open left half complex plane. 
\end{theo}

	According to Lemma 3, we need to ensure that $M'$ and $L_1 \bar{D} L_3$ are Hurwitz stable.
	The stability of matrix $M'$ is closely related to the eigenvalues of the matrix $M'$. 
	According to the Gerschgorlin circle theorem of block matrices \cite{35},
	all eigenvalues of the matrix are also located in the left half complex plane as long as all Gerschgorlin circles (15) are placed in the left half complex plane.
	Therefore, for each agent $i$, we can design $K_ i$ to adjust the position of the corresponding Gerschgorlin circle so that it is located in the open left half complex plane. 
	That is, if the control gain matrices $K_ i. i=1,2, \cdots, N $ are selected to make such that Gerschgorlin circles (15) are all located in the open left half complex plane, matrix $M'$ is Hurwitz.
	If matrix $L_1 \bar{D} L_3$ is also Hurwitz at this time, then the leader-following state consensus of MAS (1) and (2) under control protocol (12) is achieved.

\begin{rmk}
	The information used in designing the feedback gain matrix $K_i$ includes the coefficient matrix information of agent $i$ and the neighbor information it communicates with, but without using other agent information that is not directly communicates with it in the network.
	From this perspective, the design method proposed in Theorem 2 is still decentralized.
	This result improves the conclusion given in \cite{31} by solving matrix inequalities that contain global network information. 
\end{rmk}

As shown in the proof of Theorems 1 and 2, the conclusion given through the linear transformation method does not require whether the weight on the edge is a symmetric matrix, which is different from the assumption that the weighted matrix on the edge is positive definite or negative definite commonly used in most existing research on matrix weighted MASs \cite{27,28}. So to our knowledge, Theorem 1 and Theorem 2 provide a decentralized design method for consensus protocols in heterogeneous MASs under the condition that the weight matrix is asymmetric for the first time. 
The specific selection process of the gain matrix $K_i$ will be further shown in Example 3 in Section \uppercase\expandafter{\romannumeral4}.

\section{Numberical examples}
In this section, we provide several numerical examples to verify the effectiveness of the proposed method and the obtained results. Firstly, the following example is used to verify the decentralized design approach in Theorem 1.
\begin{exm}
	\begin{figure}
		\centerline{\includegraphics[width=0.3\linewidth]{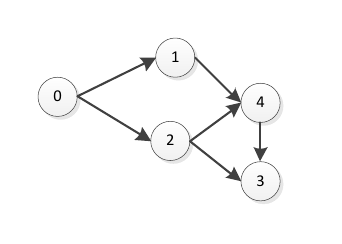}}
		\caption{Communication topology in Example 1}
		\label{fig}
	\end{figure}
	Consider the MAS with the communication topology shown in Figure 1, where node 0 is the leader and nodes from 1 to 4 are the followers.
	The parameter matrices for their dynamics are as follows
	$$
	\begin{aligned}
		& A_0=\left[\begin{array}{ll}
			2 & 0 \\
			0 & 4
		\end{array}\right],  B_0=\left[\begin{array}{l}
			1 \\
			1
		\end{array}\right] 
		A_1=\left[\begin{array}{ll}
			1 & 0 \\
			3 & 4
		\end{array}\right],  B_1=\left[\begin{array}{l}
			2 \\
			1
		\end{array}\right],  A_2=\left[\begin{array}{ll}
			1 & 2 \\
			1 & 4
		\end{array}\right],  B_2=\left[\begin{array}{l}
			2 \\
			3
		\end{array}\right]  \\
		& A_3=\left[\begin{array}{ll}
			1 & 1 \\
			0 & 4
		\end{array}\right],  B_3=\left[\begin{array}{l}
			1 \\
			4
		\end{array}\right],  A_4=\left[\begin{array}{ll}
			2 & 2 \\
			3 & 5
		\end{array}\right],  B_4=\left[\begin{array}{l}
			2 \\
			4
		\end{array}\right] \\
		&
	\end{aligned}
	$$
	The weight matrices with respect to the communication in Figure 1 are as follows
	$$
	\begin{aligned}
		& W_{10}=\left[\begin{array}{ll}
			1 & 2 \\
			3 & 4
		\end{array}\right], W_{20}=\left[\begin{array}{ll}
			4 & 3 \\
			3 & 4
		\end{array}\right], W_{41}=\left[\begin{array}{ll}
			1 & 2 \\
			5 & 4
		\end{array}\right],
		W_{42}=\left[\begin{array}{ll}
			2 & 2 \\
			3 & 4
		\end{array}\right],  W_{32}=\left[\begin{array}{ll}
			1 & 2 \\
			0 & 4
		\end{array}\right], W_{34}=\left[\begin{array}{ll}
			1 & 0 \\
			3 & 4
		\end{array}\right]
	\end{aligned}
	$$
	and $D_1=W_{10},D_2=W_{20}$.
	
	Consider the DST with edges ${(0,1)(0,2)(2,3)(2,4)}$.
	Nodes $1$ and $2$ directly receive the leader's information, and nodes $3$ and $4$ have a common neighbor node $2$ on the DST.
	
	For the design of the gain matrices, first, select $G_i$ as follows
	$$
	\begin{aligned}
		& G_1=\left[\begin{array}{ll}
			1 & 1
		\end{array}\right],\; G_2=\left[\begin{array}{ll}
			2 & 2
		\end{array}\right],
		 G_3=\left[\begin{array}{ll}
			3 & 3
		\end{array}\right],\; G_4=\left[\begin{array}{ll}
			4 & 4
		\end{array}\right]
	\end{aligned}
	$$
	such that $L_1 \bar{D} L_3$ is Hurwitz stable. 
	Then select $K_i$ according to Theorem 1
	$$
	\begin{aligned}
		& K_1=\left[\begin{array}{ll}
			0.7 & 0.1
		\end{array}\right], \; K_2=\left[\begin{array}{ll}
			-2.786 & 2.464
		\end{array}\right],
		K_3=\left[\begin{array}{ll}
			2 & -1.625
		\end{array}\right],\;  K_4=\left[\begin{array}{ll}
			-12.25 & 6
		\end{array}\right]
	\end{aligned}
	$$
	such that 
	$$A_1-B_1 G_1 -B_1 K_1 D_1,\; A_2-B_2 G_2-B_2 K_2 D_2$$
	$$A_3-B_3 G_3 -B_3 K_3 W_{32}, \; A_4-B_4 G_4 -B_4 K_4 W_{42}$$
	are Hurwitz stable.
	
	Figure 2 shows the state errors, both for the initial states 
	$ x_{0}= (2.5, 2.5), x_{1}=(0.5, 0), x_{2}= (1, 0), x_{3}=(1.5, 0), x_{4}=(2, 0)$.
	\begin{figure}
		\centerline{\includegraphics[width=0.3\linewidth]{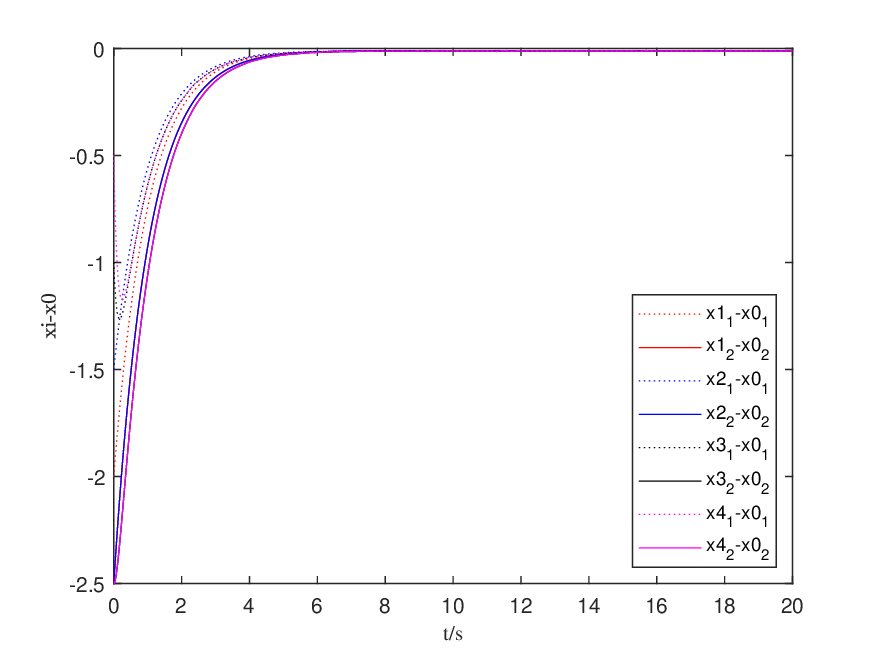}}
		\caption{State errors in Example 1}
		\label{fig}
	\end{figure}
	
	\begin{rmk}
		When solving the feedback gain matrix $K_i$, only the information of agent $i$ and the communication weight between itself and its parent node are used.
		Therefore, the solving process is decentralized.
	\end{rmk}
\end{exm}
As mentioned above, the gain matrices $G_i$ and $K_i$ play different roles in the consensus process.
Below, we will use an example to demonstrate the impact of the gain matrix $G_i$ on the consensus process.
\begin{exm}
	Consider the same MAS in Example 1.
	Now we set $G_i$  to be zero matrices.
	One can verify $L_1 \bar{D} L_3$ is Hurwitz stable.
	
	Next, we select $K_i$ according to Theorem 1
	$$
	\begin{aligned}
		& K_1=\left[\begin{array}{ll}
			0.2222 & 0.6111
		\end{array}\right], \; K_2=\left[\begin{array}{ll}
			-2.5 & 2.75
		\end{array}\right],
		 K_3=\left[\begin{array}{ll}
			2.5 & -1.10313
		\end{array}\right], \; K_4=\left[\begin{array}{ll}
			-10.25 & 6
		\end{array}\right]
	\end{aligned}
	$$
	such that
	$$
	\begin{aligned}
		&  A_1-B_1 K_1 D_1, \; A_2-B_2 K_2 D_2, A_3-B_3 K_3 W_{32},\; A_4-B_4 K_4 W_{42}
	\end{aligned}
	$$ 
	are Hurwitz stable.

	\begin{figure}
		\centerline{\includegraphics[width=0.3\linewidth]{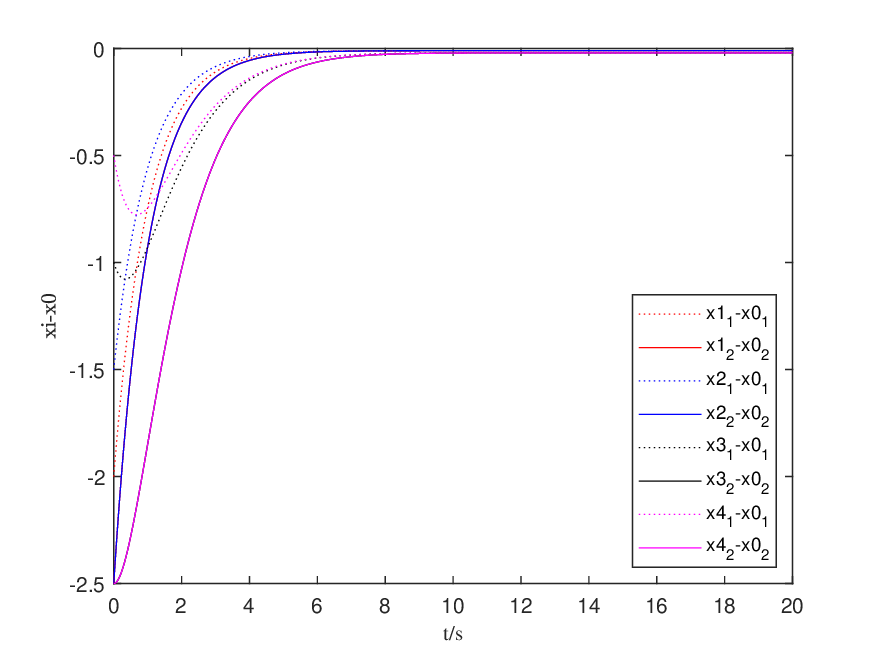}}
		\caption{State errors in Example 2}
		\label{fig}
	\end{figure}
	
	Figure 3 shows the state errors with initial states
	$ x_{0}= (2.5, 2.5), x_{1}=(0.5, 0), x_{2}=(1, 0), x_{3}=(1.5, 0), x_{4}=(2, 0)$.
	
	By comparing Example 1 and Example 2, one can observe that the gain matrices $G_i$ have an impact on the convergence rate of the consensus process.
\end{exm}

Next, we will demonstrate the decentralized design method under protocols that use all the neighbors' information.
\begin{exm}
	Consider the same MAS as in Example 1.
	We first select $G_i$
	$$
	\begin{aligned}
		& G_1=\left[\begin{array}{ll}
			1 & 1
		\end{array}\right], G_2=\left[\begin{array}{ll}
			2 & 2
		\end{array}\right],
		 G_3=\left[\begin{array}{ll}
			3 & 3
		\end{array}\right], G_4=\left[\begin{array}{ll}
			4 & 4
		\end{array}\right]
	\end{aligned}
	$$
	such that $L_1 \bar{D} L_3$ is Hurwitz stable. 
	According to Theorem 2, we design matrices $K_i$ making the following Gerschgorlin circles ($i=1,2,3,4$)
	$$
	\left(\left\|\left(A_i^{\prime}-B_i^{\prime}-\lambda I_n\right)^{-1}\right\|\right)^{-1} \leqq \sum_{j \neq 1, j \in N_i}\left\|-B_i K_i W_{i, j}\right\|
	$$
	to be all located in the open left half complex plane, where $\lambda$ are the eigenvalues of matrix $M'$.
	Thus $K_i$ can be selected  as 
	$$
	\begin{aligned}
		& K_1=\left[\begin{array}{ll}
			0.7 & 0.1
		\end{array}\right], \; K_2=\left[\begin{array}{ll}
			-2.786 & 2.464
		\end{array}\right] \\
		& K_3=\left[\begin{array}{ll}
			2 & -1.625
		\end{array}\right], \; K_4=\left[\begin{array}{ll}
			-12.25 & 5
		\end{array}\right]
	\end{aligned}
	$$
	\begin{figure}
		\centerline{\includegraphics[width=0.3\linewidth]{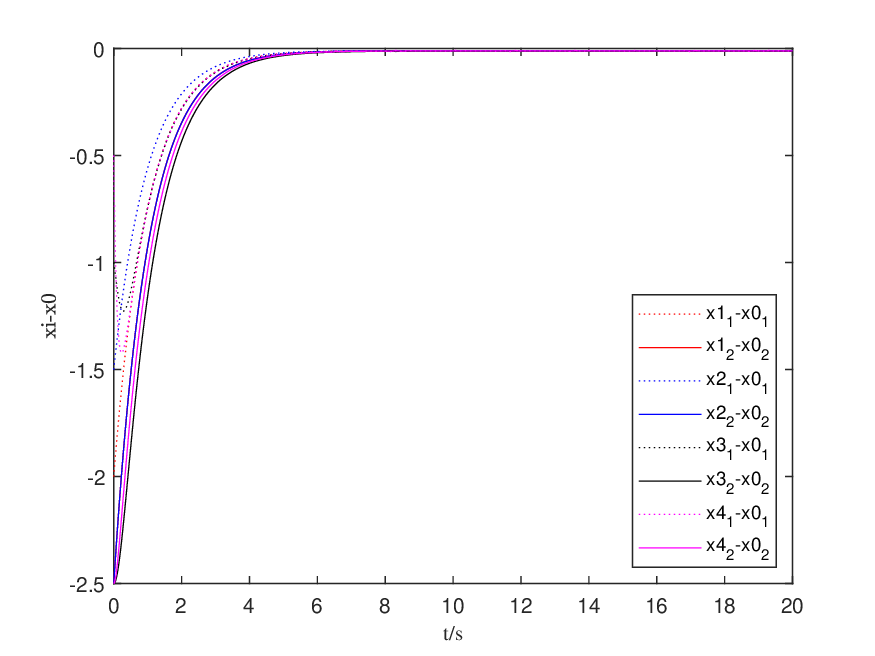}}
		\caption{State errors in Example 3}
		\label{fig}
	\end{figure} 
	Figure 4 shows the state errors with initial states
	$x_{0}=(2.5, 2.5), x_{1}=(0.5, 0), x_{2}=(1, 0), x_{3}=(1.5, 0), x_{4}=(2, 0)$.
	
	\begin{rmk}
		When solving the gain matrix $K_i$, only the information of agent $i$ and the communication weight between itself and its neighbor node are used.
		Therefore, the solving process is decentralized.
	\end{rmk}
\end{exm}

\section{CONCLUSION}
This paper studies the decentralized design of consensus protocols for heterogeneous MASs with asymmetric weight matrices.
A consensus protocol with minimal communication links is designed, and a decentralized design approach for the gain matrices in the protocol is obtained.
The decentralized design approach is further extended to the cases of protocols using all neighbor information and interconnected MASs.

The proposed approach provides a basis for dealing with more complex consensus problems of matrix-weighted MASs.
In future research, we will consider decentralized design methods for more complex situations such as switching communication topologies or involving nonlinear agents and so on.



\section*{Acknowledgement}
This work is supported by Beijing Natural Science Foundation (4232041) and National Natural Science Foundation of China (62273014).


\end{document}